\documentclass[aos,MSNbibl,nameyear,seceqn,dvips]{arximspdf}
\usepackage{algorithm,algorithmic}
\usepackage{dcolumn}
\usepackage{graphicx}

\doi{10.1214/12-AOS1025} 
\volume{40}
\issue{3}
\pubyear{2012}
\firstpage{1878}
\lastpage{1905}

\makeatletter

\newcolumntype{k}[1]{D{,}{}{#1}}

\newcommand{\textscc}[1]{\mathrm{\uppercase{#1}}}
\renewcommand{\vartheta}{\theta}

\newcommand{\cB}{\mathcal{B}}
\newcommand{\cH}{\mathcal{H}}
\newcommand{\cK}{\mathcal{K}}
\newcommand{\cN}{\mathcal{N}}
\newcommand{\cX}{\mathcal{X}}

\newcommand{\bI}{\mathbf{I}}
\newcommand{\bX}{\mathbf{X}}
\newcommand{\bY}{\mathbf{Y}}
\newcommand{\bone}{\mathbf{1}}

\newcommand{\kf}{\mathfrak{f}}

\newcommand{\R}{\mathbb{R}}
\newcommand{\p}{\mathbb{P}}
\newcommand{\E}{\mathbb{E}}

\newcommand{\ee}{\mathsf{e}}
\newcommand{\ff}{\mathsf{f}}
\newcommand{\pP}{\mathsf{P}}

\newcommand{\bxi}{\bolds{\xi}}
\newcommand{\MSE}{\operatorname{MSE}}
\newcommand{\hMSE}{\widehat{\operatorname{MSE}}}
\newcommand{\argmin}{\mathop{\operatorname{argmin}}}

\newcommand{\eps}{\varepsilon}

\newcommand{\tl}{{\tilde{\lambda}}}
\newcommand{\lth}{{\lambda_{\vartheta}}}
\newcommand{\lproj}{\lambda^{\mathrm{PROJ}}}
\newcommand{\lexp}{\lambda^{\mathrm{EXP}}}

\newtheorem{theorem}{Theorem}[section]
\newproclaim{ex}{Example}[section]
\newproclaim{exo}{Exercise}[section]
\newtheorem{prop}{Proposition}[section]
\newtheorem{cor}{Corollary}[section]
\newtheorem{lem}{Lemma}[section]
\newproclaim{rem}{Remark}[section]
\newproclaim{rems}{Remarks}[section]
\newproclaim{hyp}{Assumption}

\newproclaim{defin}{Definition}[section]
\newproclaim{notation}{Notation}

\newcommand{\eqref}[1]{(\ref{#1})}
\makeatother

\begin{document}
\begin{frontmatter}

\title{Deviation optimal learning using greedy $Q$-aggregation}
\runtitle{Deviation optimal learning}

\begin{aug}
\author[A]{\fnms{Dong} \snm{Dai}\ead[label=e1]{dongdai@stat.rutgers.edu}},
\author[B]{\fnms{Philippe} \snm{Rigollet}\corref{}\ead[label=e2]{rigollet@princeton.edu}\thanksref{t2}}
\and

\author[A]{\fnms{Tong} \snm{Zhang}\ead[label=e3]{tzhang@stat.rutgers.edu}}
\thankstext{t2}{Supported in part by NSF Grants DMS-09-06424 and
CAREER-DMS-1053987.}
\runauthor{D. Dai, P. Rigollet and T. Zhang}
\affiliation{Rutgers University, Princeton University and Rutgers University}
\address[A]{D. Dai\\
T. Zhang\\
Department of Statistics\\
Rutgers University\\
Piscataway, New Jersey 08854\\
USA\\
\printead{e1}\\
\phantom{E-mail:\ }\printead*{e3}}
\address[B]{P. Rigollet\\
Department of Operations Research \\
\quad and Financial Engineering\\
Princeton University\\
Princeton, New Jersey 08544\\
USA\\
\printead{e2}}
\end{aug}

\received{\smonth{3} \syear{2012}}
\revised{\smonth{6} \syear{2012}}

%
\begin{abstract}
Given a finite family of functions, the goal of model
selection aggregation is to
construct a procedure that mimics the function from this
family that is the closest to an unknown regression function.
More precisely, we consider a general regression model with fixed
design and measure the distance between functions by the mean squared
error at the design points.
While procedures based on exponential weights are known to solve the problem
of model selection aggregation in expectation, they are, surprisingly,
sub-optimal in deviation.
We propose a new formulation called $Q$-aggregation that addresses this
limitation;
namely, its solution leads to sharp oracle inequalities that are
optimal in a minimax sense.
Moreover, based on the new formulation, we design greedy
$Q$-aggregation procedures
that produce sparse aggregation models achieving the optimal rate.
The convergence and performance of these greedy procedures
are illustrated and compared with other standard methods on simulated examples.
\end{abstract}

%


\end{frontmatter}

\section{Introduction}
\label{SECintro}

Model selection is one of the major aspects of statistical learning
and, as such, has received considerable attention over the past
decades. More recently, the seminal works of \citet{Nem00} and \citet
{Tsy03} have introduced an idealized setup to study the properties of
model selection procedures independently of the models themselves. We
consider this so-called \textit{pure model selection aggregation} (or
simply MS aggregation) framework for the simple model of Gaussian
regression with fixed design.

Let $x_1, \ldots, x_n$ be $n$ given design points in a space $\cX$, and
let $\cH=\{f_1, \ldots,\break f_M\}$ be a given dictionary of real valued
functions on $\cX$. The goal is to estimate an unknown regression
function $\eta\dvtx  \cX\to\R$ at the design points based on observations
\[
Y_i = \eta(x_i) +\xi_i ,
\]
where $\xi_1, \ldots, \xi_n$ are i.i.d. $\cN(0, \sigma^2)$. Our main
results are actually stated for sub-Gaussian random variables, but
since most of the literature is available only for Gaussian noise, we
temporarily make this assumption to ease comparisons throughout the
\hyperref[SECintro]{Introduction}. The performance of an estimator $\hat\eta$ is measured
by its mean square error (MSE) defined by
\[
\MSE(\hat\eta)=\frac{1}{n} \sum_{i=1}^n
\bigl(\eta(x_i)-\hat\eta (x_i)\bigr)^2 .
\]
In the pure model selection aggregation framework, the goal is to build
an estimator $\hat\eta$ that mimics the function $f_j$ in the
dictionary with the smallest MSE. Formally, a good estimator $\hat\eta
$ should satisfy the following \textit{oracle inequality} in a certain
probabilistic sense:
%
\begin{equation}
\label{EQoi1} \MSE(\hat\eta) \le\min_{j=1, \ldots, M} \MSE(f_j) +
\Delta_{n,M}\bigl(\sigma^2\bigr) ,
\end{equation}
where the remainder term $\Delta_{n,M}>0$ should be as small as
possible. Note that oracle inequality \eqref{EQoi1} is a truly finite
sample result, and the remainder term should show the interplay between
the three fundamental parameters of the problem: the ``dimension'' $M$,
the sample size $n$ and the noise level $\sigma^2$. Most oracle
inequalities for model selection aggregation have been produced in
expectation [see the references in \citet{RigTsy12}] with notable
exceptions [\citet{Aud08,LecMen09,GaiLec11,DaiZha11,Rig12}] who
produced oracle inequalities that hold with high probability and to
which we will come back later.

From the early days of the pure model selection problem, it has been
established [see, e.g., \citet{Tsy03,Rig12}] that the smallest
possible order for $\Delta_{n,M}(\sigma^2)$ was $\sigma^2(\log M)/n$
for oracle inequalities in expectation, where ``smallest possible'' is
understood in the following minimax sense. There exists a dictionary
$\cH=\{f_1, \ldots, f_M\}$ such that the following holds. For any
estimator $\hat\eta$, there exists a regression function $\eta$ such that
\[
\E\MSE(\hat\eta) \ge\min_{j=1, \ldots, M} \MSE(f_j) + C
\sigma^2 \frac
{\log M}{n}
\]
for some positive constant $C$. Moreover, it follows from the same
results that this lower bound holds not only in expectation but also
with positive probability.

The established terminology \textit{model selection} is somewhat
misleading. Indeed, while the goal is to mimic the best model in the
dictionary~$\cH$, it has been shown [see \citeauthor{RigTsy12} (\citeyear{RigTsy12}),
Theorem~2.1]
that there exists a dictionary~$\cH$ such that any estimator (selector)
$\hat\eta$ restricted to be one of the elements of $\cH$ cannot
satisfy an oracle inequality such as~\eqref{EQoi1} with a~remainder
term of order smaller than $\sigma\sqrt{(\log M)/n}$, which is clearly
suboptimal. Rather than model selection, \textit{model averaging} has
been successfully employed to derive oracle inequalities in expectation
such as~\eqref{EQoi1}. More precisely, model averaging consists in
choosing $\hat\eta$ as a convex combination of the~$f_j$s with
carefully chosen weights. Let $\Lambda$ be the flat simplex of~$\R^M$
defined by
\[
\Lambda= \Biggl\{ \lambda=(\lambda_1, \ldots, \lambda_M)^\top
\in\R^M \dvtx \lambda_j\ge0 , \sum
_{j=1}^M \lambda_j=1 \Biggr\} .
\]
Each $\lambda\in\Lambda$ yields an \textit{aggregate} estimator
$\hat
\eta=\ff_\lambda$, where
\[
\ff_\lambda=\sum_{j=1}^M
\lambda_j f_j .
\]
This is why we refer to this problem as \textit{model selection
aggregation} or MS aggregation.
The early papers of \citet{Cat99} and \citet{Yan99} introduced and
proved optimal theoretical guarantees for a model averaging estimator
called \textit{progressive mixture} that was later studied in \citet
{Aud08} and \citet{JudRigTsy08} from various perspectives. This
estimator is based on \textit{exponential weights}, which, since then,
have been predominantly used and have led to optimal oracle
inequalities in expectation. Let $\pi=(\pi_1, \ldots, \pi_M)^\top
\in
\Lambda$ be a given \textit{prior} and $\beta>0$ be a temperature
parameter, then the $j$th exponential weight is given by
%
\begin{equation}
\label{EQexpweights} \lambda_j^{\textscc{exp}} \propto
\pi_j \exp \bigl(- n \hMSE (f_j)/\beta \bigr) ,
\end{equation}
where
\[
\hMSE(f_j)=\frac{1}{n}\sum_{i=1}^n
\bigl(Y_i - f_j(x_i) \bigr)^2
.
\]
The most common prior choice is the uniform prior $\pi=(1/M, \ldots,
1/M)^\top$, but other choices that put more or less weight on different
functions of the dictionary have been successfully applied to various
related problems; see, for example, \citet{DalSal11}, \citeauthor{RigTsy11} (\citeyear{RigTsy11,RigTsy12}).
Note that progressive mixture contains an extra averaging
step which is irrelevant to the fixed design problem that we study
here, but we implement it in Section~\ref{SECnum} for comparison with
the nonaveraged procedure.

The fixed design Gaussian regression was considered in Dalalyan and Tsybakov
(\citeyear{DalTsy07,DalTsy08})
who proved an oracle inequality of the form~\eqref{EQoi1}
with optimal remainder term. This result suffers two deficiencies:
first, it can be extended to other types of noise, but not to
sub-Gaussian distributions in general. Second, and perhaps most
importantly, it holds only in expectation and not with high
probability. While this second limitation may have followed the proof
technique, we actually show in Section~\ref{SECmain} that it is
inherent to exponential weights. Consequently, we say that exponential
weights are \textit{deviation suboptimal} since the expectation of the
resulting MSE is of the optimal order, but the deviations around the
expectation are not. Note also that the original paper of \citet
{DalTsy07} made some boundedness assumption on the distance between
function in the dictionary $\cH$ and the regression function~$\eta$.
This assumption was lifted in their subsequent paper [\citet
{DalTsy08}]. In this paper, we make no such assumption except for the
lower bound, which, of course, makes our result even stronger.

For regression with random design, \citet{Aud08} observed also that
various progressive mixture rules are deviation suboptimal. In the same
paper, he addressed this issue by proposing the STAR algorithm which is
optimal both in expectation and in deviations under the uniform prior
and, remarkably, does not require any parameter tuning as opposed to
progressive mixture rules. Also for random design, \citet{LecMen09}
followed by \citet{GaiLec11} proposed deviation optimal methods based
on the same sample splitting idea. However, sample splitting method do
not carry to fixed design. Subsequently, \citet{Rig12} proposed a new
estimator, similar to the one studied in the rest of the paper and that
enjoys the same theoretical properties as the STAR algorithm but for
fixed design regression. However, while it is the solution of a convex
optimization problem, Rigollet's method comes without implementation.
Finally, a~first implementation of a greedy algorithm that enjoys
optimal deviation was proposed in \citet{DaiZha11}. Our subsequent
results extend both the results of \citet{Rig12} and \citet{DaiZha11}
in various directions.

In Section~\ref{SECLB} of the present paper, we study the deviation
suboptimality of two commonly used aggregate estimators: the aggregate
by exponential weights and the aggregate by projection. Then, in
Section~\ref{SECmain}, we extend the original method of \citet{Rig12}
in several directions. First and foremost, our extension allows us to
put a prior weight on each element of the dictionary. These prior
weights appear explicitly in the oracle inequalities that are derived
in Section~\ref{SECmain}. Both the method of \citet{Rig12} and ours
are solutions of convex optimization problems. In Section~\ref
{SECalg}, we propose efficient greedy model averaging (GMA) procedures that
approximately solve the newly proposed $Q$-aggregation formulations. It
is shown that GMA can produce sparse
model aggregates that achieve optimal deviation bounds.
The performance of different model selection and aggregation estimators
are compared in Section~\ref{SECnum}.

\begin{notation*}
For any vector $v$, we denote by $v_j$ its $j$th coordinate. Moreover,
for any functions $f,g \dvtx  \cX\to\R$, we define the pseudo-norm
\[
\|f\|^2=\frac{1}{n}\sum_{i=1}^n
f(x_i)^2 ,
\]
and the associated inner product
\[
\langle f, g \rangle=\frac{1}{n}\sum_{i=1}^n
f(x_i)g(x_i) .
\]
Also, we define the function $\bY\dvtx  \cX\to\R$ to be any function such
that $\bY(x_i)=Y_i$. Observe that with the above notation, we have
\[
\hMSE(f)=\|\bY- f\|^2 ,\qquad \MSE(f)=\|\eta- f\|^2 .
\]
Finally, for any $p \ge1$, we denote by $|\cdot|_p$ the $\ell_p$ norm.
\end{notation*}

\section{Deviation suboptimality of commonly used estimators}
\label{SECLB}

It is well known [see, e.g., \citet{RigTsy12}] that the exponential
weights $\lambda^{\textscc{exp}}$ defined in~\eqref{EQexpweights} are
the solution of the following minimization problem:
%
\begin{equation}
\label{EQexplinint} \lambda^{\textscc{exp}} \in\argmin_{\lambda\in\Lambda} \Biggl\{
\sum_{j=1}^M \lambda_j
\hMSE(f_j)+ \frac{\beta}{n} \sum_{j=1}^M
\lambda_j \log \biggl( \frac{\lambda_j}{\pi_j} \biggr) \Biggr\} .
\end{equation}
It was shown in \citeauthor{DalTsy07} (\citeyear{DalTsy07,DalTsy08}) that for $\beta\ge4\sigma^2$, it holds
%
\begin{equation}
\label{EQdaltsy} \E\MSE (\ff_{\lambda^{\textscc{exp}}} ) \le\min_{j=1,
\ldots,
M} \biggl\{
\MSE(f_j) + \frac{\beta}{n}\log\bigl(\pi_j^{-1}
\bigr) \biggr\} .
\end{equation}
The proof of this result relies heavily on the fact that the oracle
inequality holds in expectation and whether the result also holds with
high probability arises as a natural question. While the paper
of~\citet
{Aud08} does not cover the fixed design Gaussian regression framework
of our paper and concerns exponential weights with an extra averaging
step, it contributed to the common belief that exponential weights
would be suboptimal in deviation. In particular, \citet{LecMen12}
derived lower bounds for the performance of exponential weights in
expectation when $\beta$ is chosen below a certain \textit{constant}
threshold in the case of regression with random design. Moreover, they
proved deviation sub-optimality of exponential weights when $\beta$ is
less than $\sqrt{n}/(\log n)$. However, these lower bounds rely heavily
on the fact that the design is random and do not extend to the fixed
design case. In particular, their construction uses $Y \equiv0$, which
is clearly an easy problem in the fixed design case. Proposition~\ref
{PROPlowexp} states precisely that exponential weights are deviation
suboptimal, if $\beta$ is chosen small enough and in particular if
$\beta$ is \textit{any} constant with respect to $M$ and $n$.

Another natural solution to solve the MS aggregation problem is to take
the vector of weights $\lambda^{\textscc{proj}}$ defined by
%
\begin{equation}
\label{EQproj} \lambda^{\textscc{proj}} \in\argmin_{\lambda\in\Lambda}\hMSE(
\ff_\lambda) .
\end{equation}
We call $\lproj$ the vector of \textit{projection weights} since the
aggregate estimator $\ff_{\lproj}$ is the projection of $\bY$ onto the
convex hull of the $f_j$s.

It has been established that this choice is \textit{near}-optimal for
the more difficult problem of \textit{convex aggregation} with fixed
design~[see \citet{JudNem00,Nem00,Rig12}] where the goal is to mimic
the best convex combination of the $f_j$s as opposed to simply
mimicking the best of them. More precisely, it follows from Theorem~3.5
in \citet{Rig12} that
\begin{eqnarray*}
\E\MSE (\ff_{\lambda^{\textscc{proj}}} ) &\le&\min_{\lambda
\in\Lambda}\MSE(\ff_\lambda) + 2
\sigma \sqrt{\frac{\log M}{n}}
\\
&\le&\min_{j=1, \ldots, M} \MSE(f_j)+ 2\sigma \sqrt{\frac{\log M}{n}}
,
\end{eqnarray*}
and a similar oracle inequality also holds with high probability. The
second inequality is very coarse, and it is therefore natural to study
whether a finer analysis of this estimator would yield an optimal
oracle inequality for the aggregate $\ff_{\lproj}$ both in expectation
and with high probability. This question was investigated by \citet
{LecMen09} who proved that $\ff_{\lambda^{\textscc{proj}}}$ cannot
satisfy an oracle inequality of the form~\eqref{EQoi1} with high
probability and with a remainder term $\Delta_{n,M}(\sigma^2)$ of order
smaller than $n^{-1/2}$. Their proof, however, heavily uses the fact
that the design is random, and we extend it to the fixed design case in
Proposition~\ref{PROPlowproj} below.

For both aggregates considered below, we use the following notation.
For each $j=1, \ldots, M$, we identify the functions $f_j$ on $\{x_1,
\ldots, x_n\}$ with a vector $\kf_j=(f_j(x_1), \ldots, f_j(x_n)) \in
\R^n$ where we systematically use the gothic font to identify such
vectors throughout the rest of the section. Moreover, for any vector of
weights $\lambda\in\R^M$, we write $\kf_\lambda=(\ff_\lambda(x_1),
\ldots, \ff_\lambda(x_n))$.

\subsection{Aggregate by exponential weights}
$\!\!\!$Consider the following dictionary~$\cH$. Assume that $M, n \ge3$. Let
$\ee^{(1)}=(1,0, \ldots, 0)^\top\in\R^n$ and $\ee^{(2)}=(0,1,0,
\ldots,\break 0)^\top\in\R^n$ be the first two vectors of the canonical basis of
$\R^n$. Moreover, let $\ee^{(3)}, \ldots, \ee^{(M)} \in\R^n$ be $M-2$
unit vectors of $\R^n$ that are orthogonal to both $\ee^{(1)}$ and
$\ee^{(2)}$. Let $ \kf_1, \ldots, \kf_M$ be such that
\[
\kf_1=\sigma\sqrt{ n} \ee^{(1)} ,\qquad \kf_2 =
\sigma(1+\sqrt n) \ee^{(2)} ,
\]
and for any $j=3,\ldots, M$, $\kf_j$ is defined by
\[
\kf_j =\kf_2+ \sigma\alpha_j
\ee^{(j)} ,
\]
where $\alpha_3, \ldots, \alpha_M \ge0$ are tuning parameters to be
chosen later. Moreover, take the regression function $\eta\equiv0$ so
that $\MSE(f_1)\le\MSE(f_j)$ for any \mbox{$j \ge2$}.
Observe that $ \|f_j\| \ge\sigma$ so that the following lower bounds
cannot be interpreted as artifacts of scaling the signal-to-noise ratio.

Assume that $M \ge4$ and $n \ge3$. We call low temperatures,
parameters $\beta>0$ such that
%
\begin{equation}
\label{EQbetalow} \beta\le\frac{2\sigma^2\sqrt{n}}{\log(8\sqrt{n})} .
\end{equation}
In particular the exponential weights employed in the literature on MS
aggregation use the low temperature $\beta=4 \sigma^2$; see, for
example, \eqref{EQdaltsy} above.
%

%

\begin{prop}
\label{PROPlowexp}
Fix $M\ge4, n\ge3$ and assume that the noise random variables $\xi_1,
\ldots, \xi_n$ are i.i.d. $\cN(0,\sigma^2)$. Let $\eta$ and $\cH$ be
defined as above. Then, the aggregate estimator $\ff_{\lambda
^{\textscc{exp}}}$ with exponential weights $\lambda^{\textscc{exp}}$ given
by~\eqref{EQexpweights} satisfies
\[
\MSE(\ff_{\lambda^{\textscc{exp}}}) \ge\min_{j=1,\ldots, M} \MSE (f_j) +
\frac{\sigma^2}{4\sqrt{n}} ,
\]
with probability at least $0.07$ at low temperatures, for any $\alpha_3, \ldots, \alpha_M \ge0$.

Moreover, if $M \ge8\sqrt{n}$ and for any $j \ge3$, we have
%
\begin{equation}
\label{EQalphaj} 2\sqrt{2\log(100M)} \le\alpha_j \le
n^{1/4} ,
\end{equation}
then, the same result holds at any temperature, with probability at
least $0.06$.
\end{prop}
\begin{pf} Note first that by homogeneity, one may assume that
$\sigma=1$.
Moreover, write for simplicity $\lambda=\lambda^{\textscc{exp}}$. If we
assume $\lambda_1 \le1/2$, we obtain
%
\begin{eqnarray}\label{EQprlow1}
| \kf_\lambda|_2^2 - |\kf_1|_2^2&
\ge&\bigl|\lambda_1 \kf_1 + (1-\lambda_1)
\kf_2\bigr|_2^2 - |\kf_1|_2^2\nonumber
\\
&=& (1-\lambda_1)^2 |\kf_2|_2^2-
\bigl(1-\lambda_1^2\bigr)|\kf_1|_2^2
\nonumber
\\[-8pt]
\\[-8pt]
\nonumber
& \geq&2(1-\lambda_1)^2\sqrt{n} + \bigl[(1-
\lambda_1)^2 -\bigl(1-\lambda_1^2
\bigr) \bigr]n
\\
&\ge&\sqrt{n}/2 -2\lambda_1 n.\nonumber
\end{eqnarray}

We first treat the low temperature case where $\beta$ is chosen as
in~\eqref{EQbetalow}. Define the event
\[
E=\bigl\{n\hMSE(f_2) + 2\sqrt{n} \le n\hMSE(f_1)\bigr\}
,
\]
and observe that $\eta\equiv0$ gives
%
\begin{equation}
\label{EQdefE} E= \bigl\{2\langle\kf_2 - \kf_1, \xi
\rangle_2 \ge| \kf_2 |_2^2 - |
\kf_1 |_2^2 + 2\sqrt{n} \bigr\} .
\end{equation}
On the one hand, we have $| \kf_2 |_2^2 - | \kf_1 |_2^2 = 1+2\sqrt{n}$,
and on the other hand
\[
| \kf_2 - \kf_1 |_2^2 =|
\kf_2|_2^2+| \kf_1
|_2^2=(2n+2\sqrt{n} +1)\ge \tfrac{1}{8}(1+4
\sqrt{n})^2 .
\]
Thus, we have
%
\begin{equation}
\label{EQboundPE} \p(E) \ge\p\bigl(2\langle\kf_2 -
\kf_1, \xi\rangle_2 \ge2\sqrt {2}|\kf_2 -
\kf_1|_2\bigr) =\p(Z\ge\sqrt2)\ge0.07 ,
\end{equation}
where $Z\sim\cN(0,1)$. In view of~\eqref{EQexpweights}, on the event
$E$, we have
\[
\lambda_1 \le\lambda_2 e^{-{2}/\beta\sqrt{n}} \le
\frac
{1}{8\sqrt
n} \le\frac{1}{2}
\]
for low temperature $\beta$ chosen as in~\eqref{EQbetalow}. Together
with~\eqref{EQprlow1}, it yields
\[
|\kf_\lambda|_2^2 - |\kf_1|^2
\ge\frac{\sqrt{n}}{4}.
\]

We now turn to the case of potentially high temperatures. Actually, the
following proof holds for \textit{any} temperature $\beta$ as long as
the $\alpha_j$s are chosen small enough.
In this case, we can expect the $M$ exponential weights to take
comparable values. To that end, define for each $j=2, \ldots, M$, the event
\[
F_j= \bigl\{\hMSE(f_j) \le\hMSE(f_1)
\bigr\} .
\]
Define $F=\bigcap_{j=2}^M F_j$, and denote by $F_j^c$ the complement of
$F_j$. Recall that $|\kf_j|_2^2= |\kf_2|_2^2 + \alpha_j^2$ so that
\begin{eqnarray*}
F_j^c &=& \bigl\{2\langle\kf_j -
\kf_1, \xi\rangle_2 \le| \kf_j
|_2^2 - | \kf_1 |_2^2
\bigr\}
\\
& = &\bigl\{2\langle\kf_2 - \kf_1, \xi
\rangle_2+ 2\langle\kf_j - \kf_2, \xi
\rangle_2 \le| \kf_2 |_2^2 - |
\kf_1 |_2^2 + \alpha_j^2
\bigr\}
\\
&\subset& E^c \cup G_j ,
\end{eqnarray*}
where the $E$ is defined in~\eqref{EQdefE}, and $G_j$ is defined as
\[
G_j= \bigl\{ 2\langle\kf_j - \kf_2, \xi
\rangle_2 \le\alpha_j^2 - 2\sqrt{n} \bigr\} .
\]
In view of~\eqref{EQalphaj}, we have
\[
\p(G_j)\le\p \bigl( 2\langle\kf_j - \kf_2,
\xi\rangle_2 \le -\alpha_j^2 \bigr) \le\p
\bigl(Z \ge\sqrt{2\log(100M)} \bigr)\le\frac
{0.01}{M} .
\]
Therefore,
\[
\p\bigl(F^c\bigr)\le\p\bigl(E^c\bigr) + \sum
_{j=2}^{M}\p(G_j) \le0.93+0.01=0.94 .
\]
Note now that on the event $F$, for any $j=2, \ldots, M$, we have
$\lambda_j \ge\lambda_1$ so that $\lambda_1 \le1/M \le1/2$.
Together with~\eqref{EQprlow1}, it yields
\[
|\kf_\lambda|_2^2 - |\kf_1|_2^2
\ge\frac{\sqrt{n}}{2}-\frac
{2n}{M} \ge \frac{\sqrt{n}}{4} ,
\]
where, in the last inequality, we used the fact that $M \ge8\sqrt{n}$.
\end{pf}

\subsection{Aggregate by projection} Our lower bound for the aggregate
by projection relies on a different construction of the dictionary. Let
$m$ be the smallest integer that satisfies $ m^2 \ge4n/13$ and let
$n,M$ be large enough to ensure that $m \ge16$, $M-1\ge2m$. Let $\ee^{(1)}, \ldots, \ee^{(m)} \in\R^n$ be the first $m$ vectors of the
canonical basis of $\R^n$. For any $j=1, \ldots, M$, the $\kf_j$s are
defined as
\[
\kf_j= \cases{
\sqrt{n}
\ee^{(j)} , & \quad $\mbox{if }  1\le j \le m,$
\vspace*{2pt}\cr
-\sqrt{n}\ee^{(j)}, &\quad  $\mbox{if } m+1\le j \le2m,$
\vspace*{2pt}\cr
0,& \quad $\mbox{if }  j=2m+1 ,$
\vspace*{2pt}\cr
\kf_1, &\quad $ \mbox{if }  j>2m+1 .$}
\]
Moreover, define $\eta\equiv0$ so that $0=\MSE(\kf_{2m+1}) \le\MSE
(\kf_j)$ for all $j\le M$.

\begin{prop}
\label{PROPlowproj}
Fix $ n\ge416, M \ge\sqrt{n}$, and assume that the noise random
variables $\xi_1, \ldots, \xi_n$ are i.i.d. $\cN(0,\sigma^2)$. Let
$\eta
$ and $\cH$ be defined as above. Then the projection aggregate
estimator $\ff_{\lambda^{\textscc{proj}}}$ with weights $\lambda^{\textscc{proj}}$ defined in~\eqref{EQproj} is such that
\[
\MSE(\ff_{\lambda^{\textscc{proj}}}) \ge\min_{j=1,\ldots, M} \MSE(f_j) +
\frac{\sigma^2}{ \sqrt{48n}} ,
\]
with probability larger than $1/4$. Moreover, the above lower bound
holds with arbitrary large probability if $n$ is chosen large enough.
\end{prop}

\begin{pf} Note first that by homogeneity, one may assume that
$\sigma
=1$.\vspace*{1pt} Next, observe that $\kf_{\lproj} =(\pP_m \xi, 0, \ldots,
0)^\top
\in\R^n$, where $\pP_m \xi\in\R^m$ is the projection of $
\tilde\xi=(\xi_1, \ldots, \xi_{m})^\top$ onto $\cB_1^m(\sqrt {n})$, the
$\ell_1$-ball of $\R^m$ with radius~$\sqrt{n}$.

Let $E$ denote the event on which $|\tilde\xi|_1 \le\sqrt{n}$ and
observe that, on this event, we have $\pP_m \xi=\tilde\xi$. It yields
\[
n\MSE(\kf_{\lproj})=\sum_{j=1}^m
\xi_j^2= |\tilde\xi|_2^2 .
\]
Let now $F$ denote the event on which $|\tilde\xi|_2^2 \ge m/2$, and
note that on $E\cap F$, it holds
\[
\MSE(\kf_{\lproj}) \ge\frac{m}{2n} \ge\sqrt{\frac{1}{13n}}.
\]

To conclude our proof, it remains to bound from below the probability
of $E\cap F$. The bounds below follow from the fact that $|\tilde\xi
|_2^2 $ follows a chi-squared distribution with $m$ degrees of freedom.
We begin by the event $E$. Using H\"older's inequality, we have
\[
\p\bigl(E^c\bigr)\le\p \biggl(|\tilde\xi|_2^2
\ge\frac{n}{m} \biggr) = \p \biggl(|\tilde\xi|_2^2
- \E|\tilde\xi|_2^2 \ge\frac{n}{m} - m \biggr).
\]
Next, using the fact that $m^2\le8n/13$ together with \citeauthor{LauMas00} [(\citeyear{LauMas00}), Lemma~1], we get
\[
\p\bigl(E^c\bigr)\le\p \biggl(|\tilde\xi|_2^2
- \E|\tilde\xi|_2^2 \ge\frac
{5m}{8} \biggr) \le
e^{-m/16} .
\]
Moreover, using~\citeauthor{LauMas00} [(\citeyear{LauMas00}), Lemma~1], we also get that
\[
\p\bigl(F^c\bigr)=\p \biggl(|\tilde\xi|_2^2
-\E|\tilde\xi|_2^2 \le-\frac
{m}{2} \biggr)\le
e^{-m/16} .
\]
Therefore, since $n \ge416$ implies $m \ge16$, we get
\[
\p(E\cap F)\ge1-\p\bigl(E^c\bigr)-\p\bigl(F^c\bigr)
\ge1-2e^{-m/16} \ge1-2/e\ge1/4 .
\]
\upqed\end{pf}

Note that we employed a different dictionary for each of the
aggregates. Therefore, it may be the case that choosing the right
aggregate for the right dictionary gives the correct deviation bounds.
In the next section, we propose a new aggregate estimator called
$Q$-aggregate, that automatically adjusts the aggregate to the
dictionary at hand.
%

\section{Deviation optimal model selection by $Q$-aggregation}
\label{SECmain}
According\break to~\eqref{EQexplinint}, the weight vector $\lexp$ considered
in the previous section minimizes a penalized linear interpolation of
the function $\lambda\to\hMSE(\ff_\lambda)$. The major novelty of the
method introduced in \citet{Rig12} compared to exponential weighting is
to add a quadratic term to this linear interpolation. We introduce a
family of estimators that extends the original estimator of Rigollet in
two directions: (i) it allows for a prior weighting of the functions in
the dictionary, and (ii) it allows to put different weight of each of
the component of the fitting criterion via the tuning parameter $\nu$
introduced below.

Let $\pi\in\Lambda$ be a given prior, and define the following
entropic penalty:
\[
\cK_\rho(\lambda, \pi)=\sum_{j=1}^M
\lambda_j \log \biggl(\frac
{\rho
(\lambda_j)}{\pi_j} \biggr) ,
\]
where $\rho$ is a real valued function on $[0,1]$ that satisfies $\rho
(t)\ge t$ such that $t\mapsto t\log\rho(t)$ is convex. We are
particularly interested in the choices $\rho= \bone$, the constant
function equal to $1$, which leads to a penalty that is linear
in~$\Lambda$, and $\rho(t)=t$, the identity function of $[0,1]$, which
leads to the well-known Kullback--Leibler penalty employed in
exponential weights.

Given a dictionary $\cH$ and observations $Y_1, \ldots, Y_n$, let $Q \dvtx
\Lambda\to\R$ be the function defined by
%
\begin{equation}
\label{EQdefQ} Q(\lambda)=(1-\nu)\hMSE(\ff_\lambda) + \nu\sum
_{j=1}^M \lambda_j \hMSE
(f_j) + \frac{\beta}{n}\cK_\rho(\lambda, \pi) ,
\end{equation}
where $\nu\in[0,1]$.
Let $\tilde\lambda\in\Lambda$ be such that
%
\begin{equation}
\label{EQdeftl} \tilde\lambda\in\argmin_{\lambda\in\Lambda} Q(\lambda) .
\end{equation}
We call $\ff_{\tilde\lambda}$ the \textit{$Q$-aggregate} estimator.
Note that on the one hand, if $\nu=1$ and $\rho(t)=t$, then $\tilde
\lambda=\lexp$, the exponential weights defined in~\eqref
{EQexpweights}. On the other hand, choosing $\nu=0$, $\rho(t)=1$ and
$\pi$ to be the uniform prior yields $\tilde\lambda=\lproj$, the
projection weights.

The next theorem shows that the $Q$-aggregate estimator is optimal both
in expectation and in deviation. It holds under less restrictive
conditions on the noise random variable $\xi_1, \ldots, \xi_n$. We say
that the random vector $\xi=(\xi_1, \ldots, \xi_n)^\top$ is
sub-Gaussian with variance proxy $\sigma^2>0$, if for all $t \in\R^n$,
its moment generating function satisfies
\[
\E\bigl[e^{t^\top\xi}\bigr]\le e^{{(\sigma^2 |t|_2^2)}/{2}} .
\]
Note that if $\xi\sim\cN_n(0, \Sigma)$, then $\xi$ is sub-Gaussian
with variance proxy given by $\sigma^2=\| \Sigma\|_{\mathrm{op}}$,
where $\|\Sigma\|_{\mathrm{op}}$ denotes the largest eigenvalue of the
covariance matrix~$\Sigma$.

Let $P$ be defined on the simplex $\Lambda$ by
\[
P(\lambda)= (1-\nu)\MSE(\ff_\lambda) + \nu\sum
_{j=1}^M \lambda_j \MSE
(f_j) .
\]

\begin{theorem}
\label{THmainUB}
Fix $\nu\in(0,1)$ and $\pi\in\Lambda$. Moreover, assume that the
noise random variables $\xi_1, \ldots, \xi_n$ are independent and
sub-Gaussian with variance proxy $\sigma^2$. Then for any $\beta\ge
\frac{2\sigma^2}{\min(\nu, 1-\nu)}$ and any $\delta\in(0,1)$, the
$Q$-aggregate estimator $\ff_{\tilde\lambda}$ satisfies
\[
\MSE(\ff_{\tilde\lambda}) \le\min_{\lambda\in\Lambda} \biggl\{P( \lambda) +
\frac{\beta}{n}\cK_\rho(\lambda, \pi)+ \frac{\beta}{n} \log(1/
\delta) \biggr\} ,
\]
with probability $1-\delta$. Moreover,
\[
\E\MSE(\ff_{\tilde\lambda}) \le\min_{\lambda\in\Lambda} \biggl\{P( \lambda) +
\frac{\beta}{n}\cK_\rho(\lambda, \pi) \biggr\} .
\]
\end{theorem}
Theorem~\ref{THmainUB} follows directly from Theorem~\ref{THapprox}
below, so we prove only the latter in Appendix~\ref{SECprTHmain}.

Our theorem implies that the $Q$-aggregate can compete with an
arbitrary~$\ff_\lambda$ in the convex hull
with $\lambda\in\Lambda$.
However, we are mainly interested in MS aggregation, where $\lambda$ is
at a vertex of the simplex $\Lambda$.
With $\nu\in(0,1)$, the theorem implies that the $Q$-aggregate
estimator is deviation optimal, unlike the aggregate with exponential
weights. This is explicitly stated in the following corollary, which
shows that our estimator solves optimally the problem of MS
aggregation. Its proof follows by simply restricting the infimum over
$\Lambda$ to the minimum over its vertices in Theorem~\ref{THmainUB}.\vadjust{\goodbreak}
Nonetheless, it is worth pointing out that our analysis focuses on
deviation bounds, and it does not allow us to recover~\eqref{EQdaltsy}
for the aggregate with exponential weights when $\nu=1$.

\begin{cor}
Under the assumptions of Theorem~\ref{THmainUB}, the $Q$-aggregate
estimator $\ff_{\tilde\lambda}$ satisfies
\[
\MSE(\ff_{\tilde\lambda}) \le\min_{j} \biggl\{\MSE(f_j) +
\frac
{\beta
}{n}\log \biggl(\frac{\rho(1)}{\pi_j \delta} \biggr) \biggr\} ,
\]
with probability $1-\delta$. Moreover,
\[
\E\MSE(\ff_{\tilde\lambda}) \le\min_{j} \biggl\{\MSE(f_j)
+ \frac{\beta
}{n}\log \biggl(\frac{\rho(1)}{\pi_j} \biggr) \biggr\} .
\]
\end{cor}

\begin{rem}
If we set $\rho(t)=1$ and employ the uniform prior $\pi_j=1/M, j=1,
\ldots, M$, then the optimization of the criterion $Q$ is independent
of $\beta$. In this case, we may simply set $\nu=1/2$, and the
$Q$-aggregate estimator becomes parameter free, and we recover the
original aggregate of \citet{Rig12}.
\end{rem}

\section{Algorithms}
\label{SECalg}

In the previous section, we introduced and analyzed the $Q$-aggregate
estimator. It can be easily seen that if $M$ is moderate, then it can
be computed efficiently since it requires solving the convex
optimization problem~\eqref{EQdeftl}.
The purpose of this section is to propose \textit{greedy model
averaging} (GMA) procedures that can approximately
solve the $Q$-aggregation formulation~\eqref{EQdeftl}. Moreover, GMA
leads to sparse estimators
(i.e., the resulting estimators only aggregate a small number of
dictionary functions)
that achieve the optimal deviation bounds.
These algorithms are thus appealing for their simplicity and
statistical interpretability.

\subsection{Approximate $Q$-aggregation}

Most numerical optimization algorithms do not find the exact minimum of
the objective function $Q$, but only approximate solutions. We
introduce two algorithms that minimize $Q$ approximately, with a very
specific error term for the optimization task. It relies on the
following quantity. Given a dictionary $\cH$, for any $\lambda\in
\Lambda$, let $V(\lambda)$ denote its \textit{variance on $\cH$}
and be
defined by
\[
V(\lambda)=\sum_{j=1}^M
\lambda_j \|f_j - \ff_\lambda\|^2 .
\]
For given $\eps_V, \eps>0$, we call $\ff_{\tilde\lambda_\eps}$ an
$(\eps_V, \eps)$-approximate $Q$-aggregate if the vector of weights
$\tilde\lambda_\eps\in\Lambda$ is such that
%
\begin{equation}
\label{EQapprox} Q(\tilde\lambda_\eps) \le\min_{\lambda\in\Lambda} \bigl\{
Q(\lambda) + \eps_V V(\lambda)+ \eps \bigr\} .
\end{equation}
Before going into the detailed description of the algorithms we state
a~generalization of Theorem~\ref{THmainUB} that is valid not only for
exact minimizers of $Q$ but also for approximate minimizers. Hereafter,
we use the convention $0/0=0$.\vspace*{-2pt}

\begin{theorem}
\label{THapprox}
Let $\eps, \eps_V, \nu>0$ be such that $\nu+\eps_V <1$ and fix
$\pi\in\Lambda$. Moreover, assume that the noise random variables $
\xi_1, \ldots, \xi_n$ are independent sub-Gaussians with variance
proxy $\sigma^2$.
Fix any $\theta\in(\eps_V/(\nu+\eps_V),1]$, and choose $\beta>0$
such that
%
\begin{equation}
\label{EQbetage} \beta\ge2\sigma^2 \max \biggl\{\frac{1}{ \nu-\eps_V(1-\theta
)/\theta};
\frac{1}{(1-\theta)(1-\nu-\eps_V)} \biggr\} . 
\\
\end{equation}
Then for any $\delta\in(0,1)$, any $(\eps_V, \eps)$-approximate
$Q$-aggregate estimator $\ff_{\tilde\lambda_\eps}$ satisfies
%
\begin{equation}
\label{EQOIP}
\qquad\MSE(\ff_{\tilde\lambda_\eps}) \le\min_{\lambda\in\Lambda
} \biggl\{P(
\lambda) + \eps_V V(\lambda)+ \frac{\eps}{\theta}+ \frac{\beta
}{n}
\cK_
\rho(\lambda, \pi) \biggr\} + \frac{\beta}{n}\log(1/\delta) ,
\end{equation}
with probability $1-\delta$. Moreover,
%
\begin{equation}
\label{EQOIE} \E\MSE(\ff_{\tilde\lambda_\eps}) \le\min_{\lambda\in\Lambda
} \biggl\{P(
\lambda) + \eps_V V(\lambda)+ \frac{\eps}{\theta} + \frac
{\beta}{n}
\cK_\rho(\lambda, \pi) \biggr\} .\vspace*{-2pt}
\end{equation}
\end{theorem}

\begin{rem}
If $\eps_V=0$, then~\eqref{EQbetage} reduces to $\beta\ge2\sigma^2/\min(\nu; (1-\theta)(1-\nu))$. Thus if $\nu< 1/2$, we can take
$\theta=1-\nu(1-\nu)$ and $\beta\ge2\sigma^2/\nu$. If $\nu\ge1/2$,
then for any $\theta\in(0, 1]$, we have $\min(\nu; (1-\theta
)(1-\nu
))= (1-\theta)(1-\nu)$. Furthermore, if $\eps=0$, then in the case
$\nu
\ge1/2$ we can let $\theta\to0$ and obtain Theorem~\ref{THmainUB}.\vspace*{-2pt}

\end{rem}

\begin{rem}
Theorem~\ref{THapprox} is related to PAC-Bayes-type inequalities that
also employ entropy regularization.
In particular, the proof involves an interpolated risk with variance
correction, and such techniques
have also appeared in earlier papers such as \citet{Audibert04} under
different context.\vspace*{-2pt}
\end{rem}

Clearly the smaller the $\eps_V$ and $\eps$, the better the oracle
inequality. Nevertheless, in the canonical example where $\pi$ is the
uniform prior, it is sufficient to have $\eps$ uniformly bounded by
$C(\log M)/n$ for some $C>0$ in order to maintain a statistical
accuracy of the same order as that of the true $Q$-aggregate.
However, if an estimator has error term with $\eps=0$ and a constant
$\eps_V>0$, then
it achieves a statistical accuracy of the same order as that of the
true $Q$-aggregate because
the variance term $\eps_V V$ vanishes at the vertices of the simplex
$\Lambda$.
This is the main reason to differentiate $\eps_V$ and $\eps$ in
\eqref
{EQapprox}. As we will show later on,
specially designed greedy algorithms can lead to an error term with
$\eps=0$, and thus such greedy algorithms achieve
optimal deviation bounds for MS aggregation.\vspace*{-2pt}

\subsection{Greedy $Q$-aggregation}

Optimizing convex functions over convex sets is the bread and butter of
modern statistical computing, with many algorithms ranging from
gradient descent to interior point (IP) methods\vadjust{\goodbreak} [see, e.g., \citet
{BoyVan04} for a recent overview]. For simple constraints sets such as
the simplex $\Lambda$ considered here, so-called \textit{proximal
methods} [see, e.g., \citet{BecTeb09}] have shown very promising
performance, especially when $M$ becomes large. However, the most
efficient of these methods (IP and proximal methods) does not output a
sparse solution in a general case.

In the sequel, we focus on simple greedy model averaging algorithms
(i.e., each iteration takes the form of a greedy
selection of a function in the dictionary)
that enjoy the following property. After $k$ iterations, these
algorithms return a vector $\lambda^{(k)}$ such that (i), $\lambda^{(k)}$
has at most $k$ nonzero coefficients, and (ii) $\ff_{\lambda
^{(k)}}$ is an approximate $Q$-aggregate estimator, where the quality
of the approximation will be made explicit. Specifically, appropriately
designed greedy algorithms can give
$\eps=0$ in~\eqref{EQapprox} for all $k \geq2$, and thus achieve
optimal deviation bounds using only
$k \geq2$ dictionary functions.

Minimizing a quadratic objective over the simplex $\Lambda$ is a common
problem in statistics and optimization.
We focus on greedy algorithms introduced into the
statistical literature by \citet{Jon92}.
In optimization, greedy algorithms over simplex $\Lambda$ are known as
\textit{Frank--Wolfe}-type (or reduced gradient) methods. Their name
refers to the original paper of \citet{FraWol56}.

\begin{algorithm}[t]
\caption{Greedy model averaging (GMA-1 and GMA-1$_{+}$)}
\label{ALGgreedy1}
\begin{algorithmic}
\REQUIRE Noisy observation $\bY$, dictionary $\cH=\{f_1, \ldots,
f_M\}
$, prior $\pi\in\Lambda$, parameters $\nu, \beta$.
\ENSURE Aggregate estimator $\ff_{\lambda^{(k)}}$.
\STATE Let $ \lambda^{(0)}=0$, $\ff_{\lambda^{(0)}}=0$.
\FOR{$k=1, 2, \ldots$}
\STATE Set $\alpha_k=\frac{2}{k+1}$
\STATE$J^{(k)}=\arg\min_j  (\nabla Q(\lambda^{(k)}) )_j$
\STATE option-1 (GMA-1) $\lambda^{(k)}= \lambda^{(k-1)}+ \alpha_k(\ee^{(J^{(k)})} - \lambda^{(k-1)})$
\STATE option-2 (GMA-1$_{+}$) $\lambda^{(k)} = \arg\min_{\lambda
\in
\Lambda} Q(\lambda) \mbox{ s.t. } \lambda_j = 0$ for $j \notin\{
J^{(1)},\ldots,J^{(k)}\}$
\ENDFOR
\end{algorithmic}
\end{algorithm}

We consider a few variants of greedy algorithms described in
Algorithms~\ref{ALGgreedy1}
and~\ref{ALGgreedy0}.
In these algorithms, $\ee^{(j)}$ denotes the $j$th vector of the
canonical basis of $\R^M$.
Both algorithms can be seen as greedy algorithms that add at most one
function from the dictionary at each iteration.
This feature is attractive as it outputs a $k$-sparse solution that
depends on at most $k$ functions from the dictionary
after $k$ iterations.
Each algorithm contains two variants: GMA-0 and GMA-0$_{+}$ in
Algorithm~\ref{ALGgreedy0},
and GMA-1 and GMA-1$_{+}$ in Algorithm~\ref{ALGgreedy1}. At the
same sparsity level~$k$,
the GMA-0$_{+}$ (resp., GMA-1$_{+}$) variant can further reduce
approximation error of GMA-0 (resp., GMA-1) in~\eqref{EQapprox} via
a\vadjust{\goodbreak}
more aggressive optimization step.
This kind of additional optimization is referred to as \textit
{fully-corrective} step [\citet{ShSrZh09}],
which is known to improve performance in practice.
The difference between Algorithms~\ref{ALGgreedy1} and~\ref
{ALGgreedy0} is that
the former uses first order information, namely the gradient $\nabla
Q$, to pick the best coordinate $J^{(k)}$ (which is the standard
Frank--Wolfe procedure in the greedy
algorithm literature), while the latter uses only zero order
information, namely, the coordinate that minimizes the objective value
$Q(\cdot)$
(which is relatively uncommon in the greedy algorithm literature). A
similar algorithm with the purpose of solving MS aggregation has
appeared in \citet{DaiZha11}.

Note that both algorithms give approximate solutions $\lambda^{(k)}$
that converges to the optimal
solution of~\eqref{EQdeftl}; that is, when $k \to\infty$, we have
$\eps_V \to0$ and $\eps\to0$ in~\eqref{EQapprox}.
The classical Frank--Wolfe style analysis of greedy algorithms leads to
the same convergence rate for both
approaches with error term of $\eps_V=0$ and $\eps>0$ in~\eqref
{EQapprox}. The result is
presented below in Proposition~\ref{PROPgma}.
Moreover, we present a new analysis that differentiates these two
algorithms. Specifically we obtain
a convergence result in Theorem~\ref{THalg0} below
with error term of $\eps=0$ in \eqref{EQapprox} for Algorithm~\ref
{ALGgreedy0}
when $k \geq2$ (but we are unable to prove the same result for
Algorithm~\ref{ALGgreedy1}).
The importance of achieving error with $\eps=0$ is that for $k \geq2$,
Algorithm~\ref{ALGgreedy0}
can produce a $k$-sparse approximate solution $\lambda^{(k)}$
of~\eqref
{EQdeftl}
that achieves optimal deviation.

\begin{algorithm}[t]
\caption{Greedy model averaging (GMA-0 and GMA-0$_{+}$)}
\label{ALGgreedy0}
\begin{algorithmic}
\REQUIRE Noisy observation $\bY$, dictionary $\cH=\{f_1, \ldots,
f_M\}
$, prior $\pi\in\Lambda$, parameters $\nu, \beta$.
\ENSURE Aggregate estimator $\ff_{\lambda^{(k)}}$.
\STATE Let $ \lambda^{(0)}=0$, $\ff_{\lambda^{(0)}}=0$.
\FOR{$k=1, 2, \ldots$}
\STATE Set $\alpha_k=\frac{2}{k+1}$
\STATE$J^{(k)}=\arg\min_j Q(\lambda^{(k-1)}+ \alpha_k(\ee^{(j)} -
\lambda^{(k-1)}))$
\STATE option-1 (GMA-0) $\lambda^{(k)}= \lambda^{(k-1)}+ \alpha_k(\ee^{(J^{(k)})} - \lambda^{(k-1)})$
\STATE option-2 (GMA-0$_{+}$) $\lambda^{(k)} = \arg\min_{\lambda
\in
\Lambda} Q(\lambda) \mbox{ s.t. } \lambda_j = 0$ for $j \notin\{
J^{(1)},\ldots,J^{(k)}\}$
\ENDFOR
\end{algorithmic}
\end{algorithm}

The following proposition follows from the standard analysis in \citet
{FraWol56,Jon92}.
It shows that the estimators $\lambda^{(k)}$ from Algorithms~\ref{ALGgreedy1}
and~\ref{ALGgreedy0} converge to the optimal solution of the
$Q$-aggregation formulation \eqref{EQdeftl}.
Therefore when $k \to\infty$, $\lambda^{(k)}$ achieves optimal
deviation bound. However,
a disadvantage of the bound is that the result does not imply optimal
deviation bounds for $\lambda^{(k)}$ when
$k$ is small (e.g., when $k=2$).
%
\begin{prop}
\label{PROPgma}
Assume that the dictionary $\cH$ is such that $\max_j \|f_j\| \le L$.
Fix $\nu\in(0,1/2)$ and $\pi\in\Lambda$. Moreover, assume that the
noise random variables $\xi_1, \ldots, \xi_n$ are independent and
sub-Gaussian with variance proxy $\sigma^2$. Take $\rho=\bone$ and
\[
\beta\ge\frac{2\sigma^2}{\nu} .
\]
Then, for any $k\ge1$, the aggregate estimator $\ff_{\lambda^{(k)}}$
where $\lambda^{(k)}$ is output by GMA-1 or GMA-0
(or GMA-1$_{+}$ or GMA-0$_{+}$) after $k$ steps, satisfies
\[
\MSE(\ff_{\lambda^{(k)}}) \le\min_{j} \biggl\{\MSE(f_j) +
\frac
{\beta
}{n}\log \biggl(\frac{1}{\pi_j \delta} \biggr) \biggr\}+
\frac
{16(1-\nu
)^2L^2}{1-2\nu}\frac{1}{k+3} ,
\]
with probability $1-\delta$. Moreover,
\[
\E\MSE(\ff_{\lambda^{(k)}}) \le\min_{j} \biggl\{\MSE(f_j)
+ \frac
{\beta
}{n}\log \biggl(\frac{1}{\pi_j} \biggr) \biggr\}+
\frac{16(1-\nu
)^2L^2}{1-2\nu}\frac{1}{k+3} .
\]
\end{prop}

\begin{rem}
For simplicity, we consider the case of $\nu< 1/2$, although similar
bounds can be obtained with $\nu\geq1/2$.
\end{rem}

\begin{rem}
The result of Proposition~\ref{PROPgma} follows from the classical
greedy algorithm analysis in \citet{FraWol56,Jon92,Barron93}. In
particular, the result for $\lambda^{(k)}$ output by GMA-1 is well
known in the literature;
see also \citet{Clarkson08,Jag11}.
For completeness, we include the proof in Appendix~\ref
{apxproof-prop-gma} especially since the
greedy step in GMA-0 (and GMA-0$_{+}$) is relatively uncommon.
\end{rem}

\begin{rem} \label{RMKfully-corrective}
It is known that the fully-corrective variants GMA-0$_{+}$ and
GMA-1$_{+}$ generally achieve better performance
than their partially-corrective counterparts GMA-0 and GMA-1 at the
same sparsity level $k$.
Although our analysis does not show their advantages,
faster convergence rates can be obtained for fully-corrective
algorithms under additional assumptions [\citet{ShSrZh09}].
Since the issue is not essential for our paper, we only illustrate the
benefit of fully-corrective updates by experiments.
\end{rem}

\begin{rem} \label{RMKoptimsimplex}
It follows from the proof of Proposition~\ref{PROPgma} that GMA-0 can
be used to optimize the function $Q$ over the simplex $\Lambda$.
Therefore, we can use it as a subroutine for \textrm{option-2} in the
description of Algorithms~\ref{ALGgreedy0} and~\ref{ALGgreedy1}. More
precisely, the following bound holds:
\[
Q\bigl(\lambda^{(k)}\bigr) \le\min_{\lambda\in\Lambda} Q(\lambda) +
\frac
{16(1-\nu)L^2}{k+3} .
\]
\end{rem}

For the approximation error $\frac{16(1-\nu)^2L^2}{1-2\nu}\frac
{1}{k+3}$ to be of the same order as the estimation error, one may
choose $k$ such that
\[
k \ge\frac{16(1-\nu)^2L^2 n}{\beta(1-2\nu) \log(1/\pi_{\max
})}-3 ,\vadjust{\goodbreak}
\]
where $\pi_{\max} =\max_j \pi_j$. In particular, if $\pi$ is the
uniform prior, then $\ff_{\lambda^{(k)}}$ solves the problem of MS
aggregation optimally after
\[
k \ge\frac{16(1-\nu)^2L^2 n}{\beta(1-2\nu) \log(M)}-3
\]
iterations.

Note that the above theorem requires
the somewhat unpleasant assumption that the functions in the dictionary
are uniformly bounded in $\|\cdot\|$ norm. Indeed, this assumption has
not appeared so far and is therefore not natural in this problem.

More importantly, the bound only leads to optimal deviation for large~$k$ of the order
$n/\log(M)$. The cause of this unpleasant issue is that the error term is
with $\eps\neq0$ and $\eps_V=0$ in~\eqref{EQapprox}.
In order to obtain optimal deviation bound, we have to derive an error
bound of the form~\eqref{EQapprox} with either $\varepsilon=O(\log
(M)/n)$, or with
$\eps= 0$ and $\eps_V \neq0$.
In the later case, we allow~$\eps_V$ to be relatively large, which
means that we do not have to solve \eqref{EQdeftl}
accurately.
The following theorem shows that such an error bound (with $\eps=0$)
can be achieved via GMA-0 (and GMA-0$_{+}$); in addition,
this result removes the assumption on the boundedness of dictionary function.

\begin{theorem}\label{THalg0}
Fix $\nu\in(0,1), k \ge2$ and $\pi\in\Lambda$. Moreover, assume
that the noise random variables $\xi_1, \ldots, \xi_n$ are independent
and sub-Gaussian with variance proxy $\sigma^2$. Take $\rho= \bone$ and
\begin{eqnarray*}
&&\beta\ge2\sigma^2\inf_{\theta\in(0,1]}\max \biggl
\{\frac{1}{ \nu- {(4(1-\nu)(1-\theta))}/{((k+3)\theta)}};\\
&&\hspace*{110pt}\frac{1}{(1-\theta
)(1-\nu)(1-
{4}/{(k+3)})} \biggr\}.
\end{eqnarray*}
%
Then the aggregate estimator $\ff_{\lambda^{(k)}}$ where $
\lambda^{(k)}$ is output by GMA-0 (or \mbox{GMA-0$_{+}$}) after $k$ steps,
satisfies
\[
\MSE(\ff_{\lambda^{(k)}}) \le\min_{j} \biggl\{\MSE(f_j) +
\frac{
\beta}{n}\log \biggl(\frac{1}{\pi_j \delta} \biggr) \biggr\} ,
\]
with probability $1-\delta$. Moreover,
\[
\E\MSE(\ff_{\lambda^{(k)}}) \le\min_{j} \biggl\{\MSE(f_j)
+ \frac{
\beta}{n}\log \biggl(\frac{1}{\pi_j} \biggr) \biggr\} .
\]
\end{theorem}

\begin{rem}
The theorem implies deviation bounds of the optimal order for all $k
\geq2$, and the constant $\beta$ decreases
to $2\sigma^2/\min(\nu,1-\nu)$ as in Theorem~\ref{THmainUB} when $k
\to\infty$. Such results indicate that the choice of $\nu$ is not
critical and any positive constant leads to the same optimal bound.
However, we can optimize the constant by choosing $\nu=1/2$ and we use
this value in the simulations.\vadjust{\goodbreak}

Moreover, a careful inspection of the proof indicates that $\ff_{\lambda^{(k)}}$ where $\lambda^{(k)}$ is output by GMA-0 (or
GMA-0$_{+}$) after $k$ steps is a $(\eps_V, 0)$-approximate
$Q$-aggregate estimator with $\eps_V=4(1-\nu)/(k+3)$. As a result, the
condition $\nu+\eps_V<1$ of Theorem~\ref{THapprox} requires that $k
\ge2$.
\end{rem}

To get a better quantitative idea of the result, we illustrate the
particular choice $\nu=1/2$. In this case, it can be easily shown that
the optimal $\theta$ is given by $\theta^\star_k=2/(\sqrt{k+3}+2)$.
Therefore, in this case, one may take
\[
\beta\ge\frac{4\sigma^2}{1-2/\sqrt{k+3}} .
\]
In particular, for $k=2$, it is sufficient to take $\beta= 20\sigma^2/(1+2/\sqrt{5})\ge37\sigma^2$.
%
Although it achieves
the optimal rate for MS aggregation, the large constant implies that it
is still beneficial to run the algorithm for
more than two iterations. This is confirmed by our experiments.

It is worth pointing out that with flat prior,
the first stage estimator $\ff_{\lambda^{(1)}}=f_{\hat{j}}$ is simply
the empirical risk minimizer
with $\hat{j}\in\argmin_j \hMSE(f_j)$. We have already pointed out
that this estimator achieves sub-optimal
deviation bounds; therefore the requirement of $k \geq2$ in our
analysis is natural.
With $k=2$, the estimator $f_{\lambda^{(2)}}$ is related to the STAR
algorithm, which can be regarded as a two-stage
greedy algorithm that minimizes the empirical loss function instead of
the $Q$-aggregation loss investigated in this paper.
This means that we cannot directly generalize the STAR algorithm to
more than two stages since it converges
to $\ff_{\lambda^{\textscc{proj}}}$ which is known to be suboptimal for
MS aggregation.

Notice that Theorem~\ref{THalg0} has consequences on optimization
problems beyond the scope of this paper. Indeed, we constructed a
greedy algorithm for which the approximation error at each iteration is
expressed as a function (here $\eps_V V$) and not simply a constant as
usual. This construction allowed us to derive convergence rate that
achieves optimal deviation bounds
for greedy model averaging, and to avoid stringent and unnatural
conditions on the boundedness of the problem. One of the key aspects of
the function $\eps_VV$ is that it vanishes on the set of vertices. We
believe that this technique may find applications in other optimization
problems.

\section{Numerical experiments}
\label{SECnum}

Although optimal deviation bounds are obtained for greedy $Q$
aggregation with
$k \geq2$, our analysis suggests that the performance can increase
when $k$ increases (due to reduced constants).
The purpose of this section is to illustrate this behavior using
numerical examples.
We focus on the average performance of different algorithms and configurations.

We identify a function $f$ with a vector $(f(x_1),\ldots, f(x_n))^\top
\in\R^n$.
Define $f_1, \ldots, f_M$ so that the $n\times M$ design matrix $\bX
=[f_1,\ldots,f_M]$ has i.i.d. standard Gaussian entries.\vadjust{\goodbreak}
Let $I_n$ denote the identity matrix of $\R^n$, and let $\Delta\sim
\cN_n(0, I_{n})$ be a random vector. The regression function is defined
by $\eta=f_1+0.5\Delta$. Note that typically $f_1$ will be the closest
function to $\eta$ but not necessarily.
The noise vector $\xi\sim\cN_n(0,4 \bI_n)$ is drawn independently
of~$\bX$.

We define the oracle model (OM) $f_{j\star}$, where $j^\star=\argmin_{j}\MSE(f_j)$.
The performance difference between an estimator $\hat{\eta}$ and the
oracle model $f_{j\star}$ is measured by the
\textit{regret} defined as
\[
R(\hat{\eta}) = \MSE(\hat{\eta}) - \MSE(f_{j\star}) .
\]

We run GMA-0, GMA-0$_{+}$, GMA-1 and GMA-1$_{+}$ algorithms for
$k$ iterations up to $k=40$.
The temperature $\beta$ of the exponential weights (EXP) is tuned via
10-fold cross-validation.
The projection aggregation (PROJ) estimator is obtained from GMA-0 with
$\nu=0$ for 250 iterations following Remark~\ref{RMKoptimsimplex}.
The fully-corrective optimization steps in GMA-0$_{+}$ and
GMA-1{$_{+}$} are implemented using GMA-0 and GMA-1
restricted to the support $\{J^{(1)},\ldots,J^{(k)}\}$ at each step
$k$. The purpose is to achieve better performance
at the same sparsity level $k$.

Since the target is $\eta=f_1 + 0.5 \Delta$, and $f_1$ and $\Delta$
are random Gaussian
vectors, the oracle model is likely $f_1$ (but it may not be $f_1$ due
to the mis-specification vector $\Delta$).
The noise $\sigma=2$ is relatively large, which implies a~situation
where the best convex aggregation does not outperform the oracle model.
This is the scenario considered in this paper.
For simplicity, all algorithms use a flat prior $\pi_j=1/M$ for all $j$.

\begin{table}[b]
\caption{Performance comparison}
\label{tableexp}
\begin{tabular*}{200pt}{@{\extracolsep{\fill}}lcc@{}}
\hline
\textbf{STAR} & \textbf{EXP} & \textbf{PROJ} \\
\hline
$0.43 \pm0.41$ & $0.386 \pm0.47$ & $0.407 \pm0.28$ \\
\hline
\end{tabular*}
\vspace*{6pt}
\begin{tabular*}{\textwidth}{@{\extracolsep{\fill
}}lk{4.6}k{4.6}k{4.6}k{4.6}k{4.6}@{}}
\hline
\multicolumn{1}{@{}l}{$\bolds{\nu=0.5}$} & \multicolumn
{1}{c}{$\bolds
{k=1}$} &
\multicolumn{1}{c}{$\bolds{k=2}$} & \multicolumn{1}{c}{$\bolds
{k=5}$} &
\multicolumn{1}{c}{$\bolds{k=20}$} & \multicolumn{1}{c@{}}{$\bolds
{k=40}$}\\
\hline
GMA-0 & 0.508, \pm0.76 & 0.42, \pm0.53 & 0.358, \pm0.42 & 0.336, \pm
0.38 & 0.332, \pm0.37 \\
GMA-0$_{+}$ & 0.508, \pm0.76 & 0.366, \pm0.5 & 0.341, \pm0.4 &
0.336, \pm0.38 & 0.336, \pm0.38 \\
GMA-1 & 0.54, \pm0.79 & 0.683, \pm0.44 & 0.391, \pm0.38 & 0.342, \pm
0.36 & 0.334, \pm0.37 \\
GMA-1$_{+}$ & 0.54, \pm0.79 & 0.381, \pm0.46 & 0.338, \pm0.38 &
0.336, \pm0.38 & 0.336, \pm0.38 \\
\hline
\end{tabular*}
\end{table}

The experiment is performed with the parameters $n=50$, $M=200$, and
$\sigma=2$, and repeated for
500 replications.
In order to avoid cluttering, the detailed regret of different
algorithms are given
in Table~\ref{tableexp-details} in the Appendix~\ref{apxtables}. Table~\ref{tableexp}
is a simplified comparison
of commonly used estimators (EXP and PROJ as well as STAR)
with GMA-0, GMA-0$_{+}$, GMA-1 and GMA-1$_{+}$ and $\nu=0.5$.
The regret is reported using the ``$\mbox{mean} \pm\mbox{standard
deviation}$'' format.

The results in Table~\ref{tableexp} indicate that for GMA-0 (or
GMA-0$_{+}$), from $k=1$ (corresponding to MS aggregation) to
$k=2$, there is significant reduction of error.
The performance of GMA-0 (or GMA-0$_{+}$) with $k=2$ is comparable
to that of the STAR algorithm.
This is not surprising as
STAR can be regarded as the stage-2 greedy model averaging estimator
based on empirical risk minimization.
We also observe that the error keeps decreasing (but at a~slower pace)
when $k >2$, which is consistent with
Theorem~\ref{THalg0}. It means that in order to achieve good
performance, it is necessary to use more
stages than $k=2$ [although this does not change the $O(1/n)$ rate for
the regret, it can significantly reduce the constant].
It becomes better than EXP when $k$ is as small as $5$, which still
gives a relatively sparse averaged model.

\begin{figure}

\includegraphics{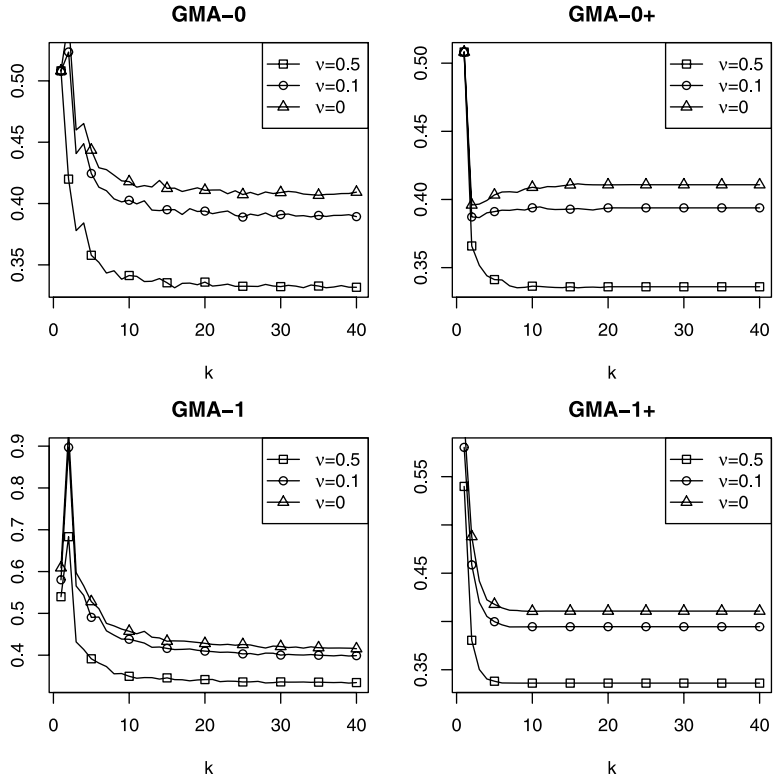}

\caption{Regrets $R(\ff_{\lambda^{(k)}})$ versus iterations $k$ for
$\nu
=0.5,0.1,0$, under 500 replications.}
\label{figexp1}
\end{figure}

Figure~\ref{figexp1} compares the MSE performance of different values
of $\nu$
for greedy algorithms considered in the paper. It shows that for the
scenario we are interested in
(i.e., where the noise is relatively large, and the best single model
is nearly as good as the best convex hull combination),
it is beneficial to choose $\nu=0.5$. Note that the greedy procedure
with $\nu=0$ converges to the convex hull
projection aggregate estimator $\ff_{\lambda^{\textscc{proj}}}$ which we
have shown to be sub-optimal for MS aggregation.
Therefore these results are consistent with our theoretical analysis,
and illustrate the importance of
$Q$-aggregation with $\nu>0$ for MS aggregation.

\begin{figure}

\includegraphics{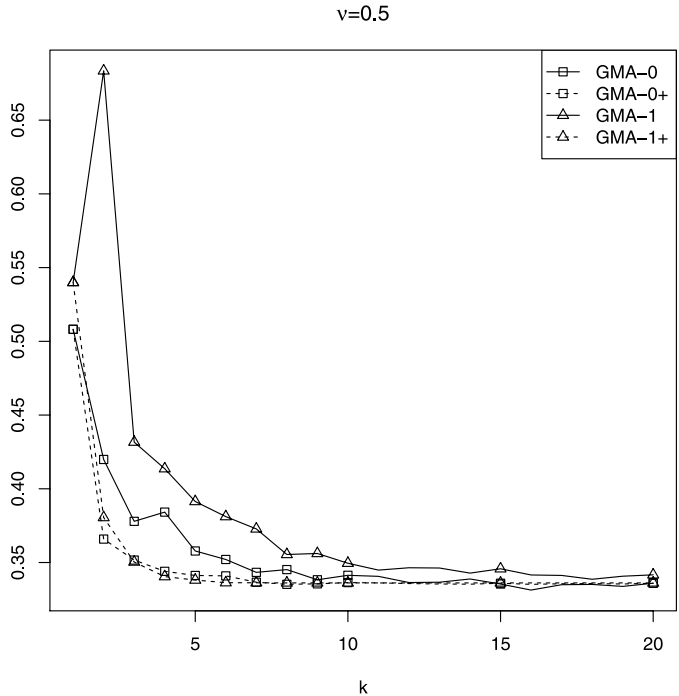}

\caption{Regrets $R(\ff_{\lambda^{(k)}})$ versus iterations $k$ of
different greedy procedures at $\nu=0.5$, under 500 replications.}
\label{figexp2}
\end{figure}

Figure~\ref{figexp2} compares the MSE of different greedy procedures
at $\nu=0.5$
(additional comparisons at $\nu=0.1$ and $\nu=0$ can be found in
Figure~\ref{figexp3} in the Appendix~\ref{apxtables}).
It shows that the classical first order greedy method \mbox{GMA-1} generally
performs worse than
GMA-0 for all $k$ and especially when~$k$ is small.
This is consistent with our theoretical analysis since Theorem~\ref
{THalg0} only applies to GMA-0.
The experiments show that the fully-corrective variants GMA-0$_{+}$
and GMA-1$_{+}$
can potentially give more accurate results than GMA-0 and GMA-1 at the
same sparsity level $k$.

\begin{appendix}\label{app}

\section{Proofs}
\label{APPA}


\subsection{\texorpdfstring{Proof of Theorem~\protect\ref{THapprox}}{Proof of Theorem 4.1}}
\label{SECprTHmain}

Let $\tilde\lambda$ be such that
\[
Q(\tilde\lambda) \le\min_{\lambda\in\Lambda} \bigl\{Q(\lambda) +
\eps_VV(\lambda) + \eps \bigr\}.
\]
Fix $\vartheta\in(0,1)$ and for any $\lambda\in\Lambda$, define
$\lth\in\Lambda$ by $\lth=(1-\vartheta)\tl+ \vartheta\lambda$.

Note that
\[
P(\tl) - P(\lth)= (1-\nu) \bigl[\MSE(\ff_\tl) -\MSE(
\ff_{\lth})\bigr] + \nu \vartheta\sum_{j=1}^M
(\tl_j-\lambda_j)\MSE(f_j). 
\]
Moreover, it is not hard to verify that
\[
\MSE(\ff_\tl)-\MSE(\ff_{\lth})=\vartheta\MSE(
\ff_\tl)- \vartheta\MSE (\ff_\lambda) + \vartheta(1-\vartheta)
\|\ff_\tl-\ff_\lambda\|^2 .
\]
The above two displays and the definition of $P(\lambda)$ yield
%
\begin{equation}
\label{EQprbase} P(\tl) - P(\lth)=\vartheta\bigl[P(\tl) - P(\lambda)\bigr] +
\vartheta (1-\vartheta ) (1-\nu) \|\ff_\tl-\ff_\lambda
\|^2 .
\end{equation}
Moreover, by the definition of $\tilde\lambda$, we have
\[
Q(\tilde\lambda) \le Q(\lth) + \eps_VV(\lth) + \eps.
\]
By replacing $Q(\tilde\lambda)$ and $Q(\lth)$ with
the expansion
\[
Q(\lambda)= P(\lambda) +\langle\bxi, \bxi\rangle- 2 \langle\bxi ,
\ff_{\lambda}-\eta\rangle + \frac{\beta}{n} \cK_\rho(\lambda,
\pi) ,
\]
where $\bxi=\bY- \eta$, we obtain
\begin{eqnarray*}
P(\tl) - P(\lth) &\le&2\langle\bxi, \ff_{\tl}-\ff_{\lth}\rangle+
\frac{\beta}{n} \cK_\rho(\lth, \pi)- \frac{\beta}{n}
\cK_\rho (\tl, \pi )+\eps_VV(\lth) + \eps
\nonumber
\\
& \le& 2\langle\bxi, \ff_{\tl}-\ff_{\lth}\rangle+
\frac{\beta
\vartheta}{n} \cK_\rho(\lambda, \pi) - \frac{\beta\vartheta}{n}
\cK_\rho(\tl, \pi)+\eps_VV(\lth)+\eps,
\end{eqnarray*}
where in the second inequality, we applied Jensen's inequality
with $\lth=(1-\vartheta)\tl+ \vartheta\lambda$
to the convex function $\lambda\mapsto\cK_\rho(\lambda, \pi)$.
Plugging \eqref{EQprbase} into this and dividing by $\vartheta$, we get
%
\begin{eqnarray}
\label{EQprA2} P(\tl) - P(\lambda) &\le&\tilde R_n(
\ff_\lambda)-(1-\vartheta ) (1-\nu) \| \ff_\tl-
\ff_\lambda\|^2
\nonumber
\\[-8pt]
\\[-8pt]
\nonumber
&&{}+ \frac{\beta}{n} \cK_\rho(\lambda,
\pi) +\frac
{\eps_V}{\vartheta}V(\lth)+\frac{\eps}{\vartheta} ,
\end{eqnarray}
where, using the fact that $\ff_\tl-\ff_{\lth}= \vartheta(\ff_\tl
-\ff_\lambda)$, we can take
\[
\tilde R_n(\ff_\lambda)=2 \langle\bxi, \ff_{\tl}-
\ff_\lambda \rangle- \frac{\beta}{n} \cK_\rho(\tl, \pi) .
\]
The following lemma allows us to control $\tilde R_n(\ff_\lambda)$ both
in expectation and with high probability.
%
\begin{lem}
\label{LEMexpDelta}
Let the noise vector $\xi=(\xi_1, \ldots, x_n)^\top$ be sub-Gaussian
with variance proxy $\sigma^2$. Then, for any $\beta>0$, $\lambda\in
\R^M$, we have
\[
\E\exp \Biggl(\frac{n }{\beta}\tilde R_n(\ff_\lambda)-
\frac
{2\sigma^2n}{\beta^2}\sum_{j=1}^M
\tl_j \Upsilon_j(\lambda) \Biggr) \le1 ,
\]
where $\Upsilon_j(\lambda)=\|f_j-\ff_{ \lambda}\|^2$.
\end{lem}
\begin{pf}
%
Fix $\lambda\in\R^M$. Using successively Jensen's inequality and the
assumption that $t \le\rho(t)$ yields
\begin{eqnarray*}
&&\E\exp \Biggl(\frac{n \tilde R_n(\ff_\lambda)}{\beta}-\frac
{2\sigma^2n}{\beta^2}\sum
_{j=1}^M \tl_j \Upsilon_j(
\lambda) \Biggr)
\\
&&\qquad= \E\exp \Biggl[ \sum_{j=1}^M
\tl_j \biggl( \frac{2 n}{\beta}\langle\bxi, f_j-
\ff_\lambda\rangle-\log \biggl(\frac
{\rho(\tl_j)}{\pi_j} \biggr) -
\frac{2\sigma^2n}{\beta^2}\Upsilon_j(\lambda) \biggr) \Biggr]
\\
&&\qquad\le\E\sum_{j=1}^M \tl_j\exp
\biggl(\frac{2 n}{\beta}\langle\bxi, f_j-\ff_\lambda\rangle-
\log \biggl(\frac{\rho(\tl_j)}{\pi_j} \biggr) -\frac{2\sigma^2n}{\beta^2}\Upsilon_j(
\lambda) \biggr)
\\
&&\qquad\le\sum_{j=1}^M \pi_j \E
\exp \biggl(\frac{2 n}{\beta}\langle \bxi, f_j-\ff_\lambda
\rangle-\frac{2\sigma^2n}{\beta^2}\Upsilon_j(\lambda ) \biggr) .
\end{eqnarray*}
Observe now that since $\xi$ is sub-Gaussian, we have
\[
\E\exp \biggl(\frac{2n}{\beta}\langle\bxi, f_j-\ff_\lambda
\rangle \biggr)\le\exp \biggl(\frac{2n\sigma^2}{\beta^2}\|\ff_\lambda-f_j
\|^2 \biggr)= \exp \biggl(\frac{2n\sigma^2}{\beta^2}\Upsilon_j(
\lambda) \biggr) .
\]
This completes the proof of our lemma.
\end{pf}

To prove the first result of Theorem~\ref{THapprox}, note that
Lemma~\ref{LEMexpDelta} together with a Chernoff bound yield that for
any $\delta\in(0,1)$,
%
\begin{equation}
\label{EQprA3} \tilde R_n(\ff_\lambda) \le\frac{2\sigma^2}{\beta}
\sum_{j=1}^M \tl_j\|
f_j-\ff_{ \lambda}\|^2+\frac{\beta\log(1/\delta)}{n}
\end{equation}
with probability at least $1-\delta$.
By combining \eqref{EQprA2} and \eqref{EQprA3}, and using the
definition of
$P(\tl)$, we obtain
%
\begin{eqnarray}\label{EQprA4}
&& (1-\nu)\MSE(\ff_\tl) + \nu\sum_{j=1}^M
\tl_j \MSE(f_j)\nonumber
\\
&&\qquad\le P(\lambda) + \frac{2\sigma^2}{\beta}\sum_{j=1}^M
\tl_j\|f_j-\ff_{ \lambda}\|^2+
\frac{\beta\log(1/\delta)}{n}
\\
&&\qquad\quad{} -(1-\vartheta) (1-\nu) \|\ff_\tl-\ff_\lambda\|^2
+ \frac{\beta
}{n} \cK_\rho(\lambda, \pi) +\frac{\eps_V}{\vartheta}V(\lth)+
\frac{\eps
}{\vartheta} .
\nonumber
\end{eqnarray}
The following identities follows directly from algebra:
\begin{eqnarray*}
\sum_{j=1}^M \tl_j
\MSE(f_j)&=&\MSE(\ff_{\tl}) + V(\tl) ,
\\
\sum_{j=1}^M \tl_j
\|f_j - \ff_\lambda\|^2&= &V(\tl) + \|
\ff_{\tl
} -\ff_\lambda\|^2 .
\end{eqnarray*}
Together with \eqref{EQprA4}, they yield
%
\begin{eqnarray}
\label{EQprA5} \MSE(\ff_\tl) & \le& P(\lambda) + \biggl(
\frac{2\sigma^2}{\beta} - \nu \biggr) V(\tl)+\frac{\beta
\log
(1/\delta)}{n}
\nonumber
\\
&&{} + \biggl[\frac{2\sigma^2}{\beta} -(1-\vartheta) (1-\nu) \biggr] \|
\ff_\tl -\ff_\lambda\|^2 + \frac{\beta}{n}
\cK_\rho(\lambda, \pi) \\
&&{}+\frac{\eps_V}{\vartheta}V(\lth)+\frac{\eps}{\vartheta}
.\nonumber
\end{eqnarray}

We now use the following identity which again follows directly from algebra:
\[
V(\lth)=\vartheta V(\lambda)+(1-\vartheta)V(\tl)+\vartheta (1-\vartheta) \|
\ff_{\tl} -\ff_\lambda\|^2 .
\]
Together with~\eqref{EQprA5}, we obtain
\[
\MSE(\ff_{\tl})\le P(\lambda)+G_1V(\tl)+G_2\|
\ff_\tl-\ff_\lambda\|^2+\frac{\beta}{n}
\cK_\rho(\lambda, \delta\pi) + \eps_V V(\lambda)+
\frac{\eps}{\vartheta} ,
\]
where $\delta\pi=(\delta\pi_1, \ldots, \delta\pi_M)^\top$,
\[
G_1=\frac{2\sigma^2}{\beta} -\nu+ \frac{\eps_V(1-\vartheta
)}{\vartheta}
\]
and
\[
G_2=\frac{2\sigma^2}{\beta} -(1-\vartheta) (1-\nu)+\eps_V
(1-\vartheta) .
\]

%

To complete the proof of~\eqref{EQOIP}, it is sufficient to note that
choosing $\beta$ as in~\eqref{EQbetage} ensures that $G_1\le0$ and
$G_2 \le0$.

Using the convexity inequality $t\le e^t-1$ for any $t\in\R$, it
yields that~\eqref{EQprA3} also holds in expectation. The proof of
\eqref{EQOIE} is then concluded in the same way as the proof of~\eqref
{EQOIP} by making statements in expectation instead of statements
that hold with high probability.

\subsection{\texorpdfstring{Proof of Theorem~\protect\ref{THalg0}}{Proof of Theorem 4.2}}

It follows from a Taylor expansion that for any $\mu, \mu' \in
\Lambda
$, we have
%
\begin{equation}
\label{EQquad} Q(\mu)= Q\bigl(\mu'\bigr) + \bigl(\mu-
\mu'\bigr)^\top\nabla Q\bigl(\mu'\bigr)+ (1-
\nu) \|\ff_{\mu} - \ff_{\mu'}\|^2 .
\end{equation}
Observe also that for any $\lambda\in\Lambda$, we have (both for
GMA-0 and GMA-0$_{+}$)
\[
Q \bigl(\lambda^{(k+1)} \bigr) \le\sum_{j=1}^M
\lambda_j Q \bigl( \lambda^{(k)}+ \alpha_{k+1}
\bigl(\ee^{(j)} - \lambda^{(k)}\bigr) \bigr) .
\]
Expanding each term on the right-hand side using \eqref{EQquad}
with $\mu=\lambda^{(k)}+ \alpha_{k+1}(\ee^{(j)} - \lambda^{(k)})$ and
$\mu'=\lambda^{(k)}$ yields
%
\begin{eqnarray}
\label{EQprapprox1} Q \bigl(\lambda^{(k+1)} \bigr) &\le& Q \bigl(
\lambda^{(k)} \bigr) +\alpha_{k+1}^2 (1-\nu)\sum
_{j=1}^M \lambda_j
\|f_j -\ff_{\lambda^{(k)}}\|^2
\nonumber
\\[-8pt]
\\[-8pt]
\nonumber
&&{}+ \alpha_{k+1}
\bigl(\lambda- \lambda^{(k)}\bigr)^\top\nabla Q \bigl(
\lambda^{(k)} \bigr) .
\end{eqnarray}
Note that
\[
\sum_{j=1}^M \lambda_j
\|f_j -\ff_{\lambda^{(k)}}\|^2=\sum
_{j=1}^M \lambda_j \|f_j -
\ff_{\lambda}\|^2 + \|\ff_{\lambda^{(k)}}-\ff_{\lambda
}
\|^2 .
\]
Moreover, applying~\eqref{EQquad} with $\mu=\lambda$ and $\mu'=
\lambda^{(k)}$ yields
\[
\alpha_{k+1}\bigl(\lambda- \lambda^{(k)}\bigr)^\top
\nabla Q \bigl( \lambda^{(k)} \bigr)=\alpha_{k+1} \bigl[ Q(
\lambda)- Q \bigl(\lambda^{(k)} \bigr) \bigr]-(1-\nu)
\alpha_{k+1}\|\ff_{ \lambda^{(k)}} -\ff_\lambda\|^2.
\]
Plugging the above two displays into~\eqref{EQprapprox1} and using
$\alpha_{k+1}^2-\alpha_{k+1}\leq0$, we get
\[
Q \bigl(\lambda^{(k+1)} \bigr) \le Q \bigl(\lambda^{(k)} \bigr)
+\alpha_{k+1}^2 (1-\nu)\sum_{j=1}^M
\lambda_j \|f_j -\ff_{\lambda}\|^2 +
\alpha_{k+1} \bigl[ Q(\lambda)- Q \bigl(\lambda^{(k)} \bigr)
\bigr] .
\]
This can be written as
%
\begin{equation}
\label{EQinduction} \delta_{k+1} \le(1-\alpha_{k+1})
\delta_k + \alpha_{k+1}^2B ,
\end{equation}
where
\[
\delta_k= Q \bigl(\lambda^{(k)} \bigr)- Q (\lambda ) ,\qquad
B=(1-\nu)\sum_{j=1}^M
\lambda_j \|f_j -\ff_{\lambda}\|^2 .
\]
To conclude that
%
\begin{equation}
\label{EQinduchyp} \delta_k \le\frac{4B}{k+3} ,
\end{equation}
we proceed by induction on $k$. It is easy to see
from \eqref{EQinduction} with $k=0$ and $\alpha_1=1$
that $\delta_1 \leq B$.

Now for $k\ge1$, bound~\eqref{EQinduction} yields
\begin{eqnarray*}
\delta_{k+1}& \le &\biggl(1-\frac{2}{2+k} \biggr)
\delta_k + \biggl(\frac
{2}{2+k} \biggr)^2B
\\
& \le& \biggl(1-\frac{2}{2+k} \biggr) \frac{4B}{k+3} + \biggl(
\frac
{2}{2+k} \biggr)^2B
 = \frac{4(k^2+3k+3)B}{(k+2)^2(k+3)}\le\frac{4B}{k+4} ,
\end{eqnarray*}
where in the second inequality, we used~\eqref{EQinduchyp}.
We have proved that for any~$\lambda$, it holds
\[
Q \bigl(\lambda^{(k)} \bigr) \le Q(\lambda) + \frac{4(1-\nu)}{k+3} \sum
_{j=1}^M \lambda_j
\|f_j -\ff_{\lambda}\|^2 .
\]

To complete the proof, we check that the assumptions of Theorem~\ref
{THapprox} with $\eps_V=4(1-\nu)/(k+3)$ and $\eps=0$ are satisfied.
Moreover, using expression~\eqref{EQbetage}, we get the desired bound
on $\beta$. To conclude, notice that $V(\lambda)$ vanishes at the
vertices of the simplex $\Lambda$.

\subsection{\texorpdfstring{Proof of Proposition~\protect\ref{PROPgma}}{Proof of Proposition 4.1}}
\label{apxproof-prop-gma}

Similarly to the proof of Theorem~\ref{THalg0}, for both GMA-1 and
GMA-1$_{+}$, we have
\begin{eqnarray*}
&& Q \bigl(\lambda^{(k+1)} \bigr)-(1-\nu)\alpha_{k+1}^2
\| f_{J^{(k)}}-\ff_{\lambda^{(k)}}\|^2
\\
&&\qquad= Q \bigl( \lambda^{(k)}\bigr) + \alpha_{k+1}\bigl(
\ee^{(J^{(k)})} - \lambda^{(k)}\bigr)^\top\nabla Q\bigl(
\lambda^{(k)}\bigr)
\\
&&\qquad\leq \sum_{j=1}^M \lambda_j
\bigl[ Q \bigl( \lambda^{(k)}\bigr) + \alpha_{k+1}\bigl(
\ee^{(j)} - \lambda^{(k)}\bigr)^\top\nabla Q\bigl(
\lambda^{(k)}\bigr) \bigr]
\\
&&\qquad= Q \bigl(\lambda^{(k)} \bigr) + \alpha_{k+1}\nabla Q \bigl(
\lambda^{(k)} \bigr)^\top\bigl(\lambda- \lambda^{(k)}
\bigr)
\\
&&\qquad= Q \bigl(\lambda^{(k)} \bigr) + \alpha_{k+1} \bigl[ Q(
\lambda )-Q \bigl(\lambda^{(k)} \bigr)-(1-\nu)\|\ff_\lambda-
\ff_{\lambda^{(k)}}\|^2 \bigr] .
\end{eqnarray*}
Using $\|f_{J^{(k)}}-\ff_{\lambda^{(k)}}\|^2 \leq4 L^2$, we obtain
%
\begin{equation}
\label{EQinduction-Bp} \delta_{k+1} \leq(1-\alpha_{k+1})
\delta_k + \alpha_{k+1}^2 B' ,
\end{equation}
where we define
\[
\delta_k= Q \bigl(\lambda^{(k)} \bigr)- Q (\lambda ) ,\qquad
B'= 4(1-\nu )L^2 .
\]
Note that \eqref{EQinduction-Bp} also holds for GMA-0 and GMA-0$_{+}$ due to \eqref{EQinduction}.
Therefore similarly to the proof of Theorem~\ref{THalg0}, we can solve
the recursion in \eqref{EQinduction-Bp}
to obtain
\[
\delta_k \le\frac{4 B'}{k+3} = \frac{16(1-\nu)L^2}{k+3} .
\]
That is, we have
\[
Q\bigl(\lambda^{(k)}\bigr) \le\min_{\lambda\in\Lambda} Q(\lambda) +
\frac
{16(1-\nu)L^2}{k+3} .
\]
We can thus apply Theorem~\ref{THapprox} with $\eps_V=0$, $\eps
=16(1-\nu)L^2/(k+3)$,
and $\theta=(1-2\nu)/(1-\nu)$
to complete the proof.

\begin{table}[b]
\tabcolsep=0pt
\caption{Regret of different algorithms: oracle model is superior to
averaged models}
\label{tableexp-details}
\begin{tabular*}{200pt}{@{\extracolsep{\fill}}lcc@{}}
\hline
\textbf{STAR} & \textbf{EXP} & \textbf{PROJ} \\
\hline
\multicolumn{1}{@{}l}{$0.43 \pm0.41$} & \multicolumn{1}{c}{$0.386
\pm
0.47$} & \multicolumn{1}{c@{}}{$0.407 \pm0.28$} \\
\hline
\end{tabular*}
\vspace*{6pt}
\begin{tabular*}{\textwidth}{@{\extracolsep{4in minus
4in}}lk{4.6}k{4.6}k{4.6}k{4.6}k{4.6}@{}}
\hline
& \multicolumn{1}{c}{$\bolds{k=1}$} & \multicolumn{1}{c}{$\bolds
{k=2}$} & \multicolumn{1}{c}{$\bolds{k=5}$} &
\multicolumn{1}{c}{$\bolds{k=20}$} & \multicolumn{1}{c@{}}{$\bolds
{k=40}$}\\
\hline
\multicolumn{1}{@{}l}{GMA-0}&&&&&\\
$\nu=0.5$ & 0.508, \pm0.76 & 0.42, \pm0.53 & 0.358, \pm0.42 & 0.336,
\pm0.38 & 0.332, \pm0.37\\
$\nu= 0.1$ & 0.508, \pm0.76 & 0.523, \pm0.5 & 0.424, \pm0.35 &
0.394, \pm0.3 & 0.389, \pm0.3\\
$\nu= 0$ & 0.508, \pm0.76 & 0.55, \pm0.48 & 0.444, \pm0.34 & 0.411,
\pm0.29 & 0.409, \pm0.28\\[3pt]
\multicolumn{1}{@{}l}{GMA-0$_{+}$}&&&&&\\
$\nu=0.5$ & 0.508, \pm0.76 & 0.366, \pm0.5 & 0.341, \pm0.4 & 0.336,
\pm0.38 & 0.336, \pm0.38\\
$\nu= 0.1$ & 0.508, \pm0.76 & 0.387, \pm0.44 & 0.391, \pm0.33 &
0.394, \pm0.3 & 0.394, \pm0.3\\
$\nu= 0$ & 0.508, \pm0.76 & 0.396, \pm0.43 & 0.403, \pm0.32 & 0.411,
\pm0.29 & 0.411, \pm0.29\\[3pt]
\multicolumn{1}{@{}l}{GMA-1}&&&&&\\
$\nu=0.5$ & 0.54, \pm0.79 & 0.683, \pm0.44 & 0.391, \pm0.38 & 0.342,
\pm0.36 & 0.334, \pm0.37\\
$\nu= 0.1$ & 0.58, \pm0.83 & 0.897, \pm0.35 & 0.49, \pm0.31 & 0.41,
\pm0.29 & 0.399, \pm0.29\\
$\nu= 0$ & 0.609, \pm0.84 & 0.937, \pm0.32 & 0.528, \pm0.3 & 0.428,
\pm0.27 & 0.415, \pm0.28\\[3pt]
\multicolumn{1}{@{}l}{GMA-1$_{+}$} &&&&&\\
$\nu=0.5$ & 0.54, \pm0.79 & 0.381, \pm0.46 & 0.338, \pm0.38 & 0.336,
\pm0.38 & 0.336, \pm0.38\\
$\nu= 0.1$ & 0.58, \pm0.83 & 0.459, \pm0.45 & 0.4, \pm0.31 & 0.395,
\pm0.3 & 0.395, \pm0.3\\
$\nu= 0$ & 0.609, \pm0.84 & 0.488, \pm0.45 & 0.418, \pm0.3 & 0.411,
\pm0.29 & 0.411, \pm0.29\\
\hline
\end{tabular*}
\end{table}

\begin{figure}[b]

\includegraphics{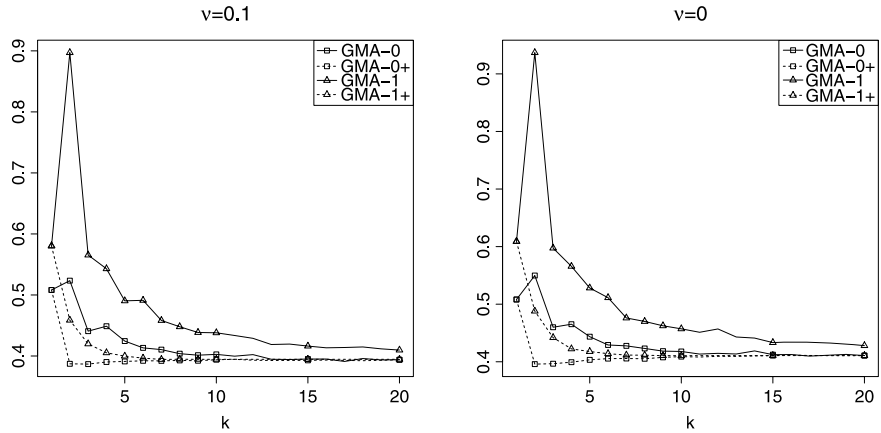}

\caption{Regrets $R(\ff_{\lambda^{(k)}})$ versus iterations $k$ of
different greedy procedures at $\nu=0.1$ and $\nu=0$, under 500 replications.}
\label{figexp3}
\end{figure}

\newpage
\section{Detailed Performance Table and Figures}
\label{apxtables}

%
%
%
\end{appendix}
\newpage

\section*{Acknowledgment}
We would like to thank an anonymous referee
for suggesting an improvement of bound~\eqref{EQbetage}.

%


\printaddresses

\end{document}